\documentclass[a4paper,10pt]{article}

\def \Z {{\mathbf {Z}}}

\def\u{\bigsqcup}
\def\eps{\varepsilon}
\title{  On Mixing  Constructions with   Algebraic Spacers }
\author{V.V. Ryzhikov\footnote{\large This work is partially  supported by  the grant
NSh 8508.2010.1.}}

\begin{document}
\Large
\maketitle

\section{Introduction.}  Following \cite{O} we  recall a proof of  the  mixing
for almost all Ornstein's  stochastic rank one constructions (section 2), then  we replace stochastic spacers by special algebraic ones \cite{R} and  in this new situation we  deduce the mixing from the weakly mixing property(section 3).
\\ 
{\bf Rank one construction} is determined by $h_1$, a sequence $r_j$  of cuttings and  a sequence $\bar s_j$ of spacers
$$ \bar s_j=(s_j(1), s_j(2),\dots, s_j(r_j-1),s_j(r_j)).$$  We recall its definition.
Let our  $T$ on the step $j$  is defined   on  a collection of disjoint sets (intervals)
$$E_j, TE_j T^2,E_j,\dots, T^{h_j-1}E_j$$
($T$ is not defined on the latest  interval $T^{h_j}E_j$).
We cut $E_j$ into $r_j$ sets (subintervals)  of the same measure
$$E_j=E_j^1\u E_j^2\u  E_j^3\u\dots\u E_j^{r_j},$$  
then for all $i=1,2,\dots, r_j$ we  consider columns
$$E_j^i, TE_j^i ,T^2 E_j^i,\dots, T^{h_j}E_j^i.$$
Adding $s_j(i)$ spacers we obtain a collection of  disjoint intervals 
$$E_j^i, TE_j^i T^2 E_j^i,\dots, T^{h_j}E_j^i, T^{h_j+1}E_j^i, T^{h_j+2}E_j^i, \dots, T^{h_j+s_j(i)}E_j^i.$$
Setting 
for all  $i<r_j$ 
$$TT^{h_j+s_j(i)}E_j^i = E_j^{i+1}$$ 
 we get  $(j+1)$-tower 
$$E_{j+1}, TE_{j+1} T^2 E_{j+1},\dots, T^{h_{j+1}}E_{j+1},$$
where 
 $$E_{j+1}= E^1_j,$$
$$h_{j+1}+1=(h_j+1)r_j +\sum_{i=1}^{r_j}s_j(i).$$
 Step by step we define a  construction $T$ on a union $X$ of all above intervals, assuming $\mu(X)=1$.

\bf On notations. \rm We denote weak operator approximations by $\approx_w$, and $\approx$ for strong ones.
$\Theta$ is the orthogonal projection into the space of constant functions in $L_2(X,\mu)$.
The  expression $T^m\approx_w\Theta$ (for  large $m$) means that $T$ is mixing. 
 
\section{Stochastic constructions}  D. Ornstein has proved \cite{O} the mixing for almost all special rank one constructions. His approach  can be  presented  in the  following manner.
Let $H_j\to \infty$,  $H_j<<r_j$.  For  uniformly distributed stochastic variables
$a_j(i)\in \{0, 1, \dots, H_j\}$ we set 
$$s_j(i)=H_j+a_j(i)- a_j(i+1).$$

Then for $m\in [h_j,h_{j+1})$
$$T^m = \hat D_1T^{m} +\hat D_2T^{m} +\hat D_3T^{m}\approx_w  \hat D_1T^{k_1}P_1 +\hat D_2T^{k_2}P_2 +\hat D_3T^{k_3}P_3,$$ 
where 

 $\hat D_i$ are operators of multiplication by indicators of  certain sets
$D_1, D_2, D_3=X\setminus (D_1\u D_2)$,
all $\hat D_i$ (and $k_i$, $P_i$) depend on $m$; 

 $k_1=m-h_{j+1}-H_{j+1},$ $|k_2|,|k_3|<h_j$;  

the operators $P_i$  \bf for almost all constructions \rm $T$  satisfy 
$$P_1\approx \sum_{n\in [-H_{j+1}, H_{j+1}]}c_{j+1}(n)T^{ n}, \ \  c_{j+1}(n)=\frac{H_{j+1}+1-|n|}{(H_{j+1}+1)^2}, $$
$$P_{2,3}\approx \sum_{n\in [-H_{j}, H_{j}]}c_{j}(n)T^{ n}, \ \  c_{j}(n)=\frac{H_{j}+1-|n|}{(H_{j}+1)^2}, $$
for all large $m$.  
As for the operators
 $P_i$, they satisfy  
$$ P_i^\ast P_i-TP_i^\ast P_i\approx_w 0,$$   
this implies  $P_i^\ast P_i\approx_w\Theta$ as $T$ is ergodic  ( all rank one  transformations are ergodic!), hence, $P_i\approx \Theta$.  From $\|\hat D_iT^{k_i}\|\leq 1$ and $P_i\approx\Theta$  we  get for a. a. constructions \rm $T$ 
$$T^m\approx_w\Theta$$
for  all large  $m$, so $T$ is mixing. 

\bf More details.  \rm
Denote
$$S_j(k,N)=-NH_j+\sum_{i=k}^{N}s_j(i)=a_j(k)-a_j(k+N).$$
For example, as   $m=Nh_j$ ( $N<<r_j$) we get
$$T^{-Nh_j}\approx_w  \frac{ 1}{r_j-N-1}\sum_{k=1}^{r_j-N-1}T^{\sum_{i=k}^{N}s_j(i)}=$$
$$=\frac{ 1}{r_j-N-1}\sum_{k=1}^{r_j-N-1}T^{NH_j+ S_j(k,N)}= \frac{ 1}{r_j-N-1}\sum_{k=1}^{r_j-N-1}T^{NH_j+a_j(k)-a_j(k+2)}.$$

If $$m=N(h_j+H_j)+v, \  N<<r_j,  \  NH_j<<h_j, \  v<<h_j,$$
then 
$$T^{m}\approx_w   \frac{ 1}{r_j-N-1}T^{k_2} \sum_{k=1}^{r_j-N-1}T^{-a_j(k)+a_j(k+2)},\eqno (1)$$
where $k_2=v-NH_j. $

Generally we define $D_1, D_2, D_3$:
$$D_1=\u_{i=0}^{m}T^i E_{j+1},$$ 
$$D_3=(X\setminus D_1) \cap\u_{i=0}^{k_2}T^i E_{j},$$
$$h_j> k_2=m-N(h_j+H_j),$$
$$D_2=X\setminus (D_1\cup D_2),$$
and use the following  approximation:  
$$T^m = \hat D_1T^{m} +\hat D_2T^{m} +\hat D_3T^{m}\approx_w  \hat D_1T^{k_1}Q_1 +\hat D_2T^{k_2}Q_2 +\hat D_3T^{k_3}Q_3,$$
where  $k_1=m-h_{j+1}-H_{j+1}$, $k_3+k_2=h_j+H_j$, and 
$$Q_1  
=\frac{1}{r_{j+1}-1}\sum_{i=1}^{r_{j+1}-1}T^{-a_{j+1}(i)+a_{j+1}(i+1)},$$
\vspace{2mm}
$$Q_2 
=\frac{1}{r_j-N}\sum_{i=0}^{r_j-N}T^{-a_j(i)+a_j(i+N)},$$
\vspace{2mm}
$$Q_3 
=\frac{1}{r_j-N-1}\sum_{i=0}^{r_j-N-1}T^{-a_j(i)+a_j(i+N+1)}$$
(for the biggest $N$  satisfied $Nh_j< m$).
 
Let's remark that sometimes   (as  in (1)) some $D_i$  becomes  of a  small measure ($\|\hat D_i\|\approx 0$).
This time the corresponding operator $Q_i$ could be out of consideration.  

\vspace{2mm}

\ \ \ \ \ \ \ \ \ \ \ \ \ \it  If $Q_i\approx \Theta$\ ($i=1,2,3)$, then $T$ is mixing. \rm
\vspace{2mm}

\bf For almost all stochastic $T$ \rm 
for a  vector $\{a_j(i)\}$,  $i \in [1, r_j]$,   the frequency 
  $$\frac{|\{i \in [1, r_j-N]  \ :\ a_j(i)-a_j(i+N)=n\}|}{r_j-N}$$ is close  to $c_j(n)$ 
(we recall that  $H_j << r_j$). Here we assume that  $N<(1-\delta_j)r_j$ for 
$\delta_j\to 0$ very slowly.  This explains the above approximations
$Q_i\approx P_i$.     
\section{ Algebraic spacers instead of  stochastic ones}
Now we present  certain effective spacer sequences and  another arguments to get  
$$Q_i\approx \Theta.$$

 Let $r_j$ be prime, $r_j\to\infty$.  We fix   generators $q_j$ in the multiplicative groups (associated with the sets $\{1, 2,\dots, r_j-1\}$) of the  fields $\Z_{r_j}$.
For some sequence  $\{H_j\}$, $H_j\geq r_j$, we define a spacer sequence
$$s_j(i)=H_j + \{q_j^i\} -\{q_j^{i+1}\},  \ \ i=1, 2, \dots,  r_j-1,$$
 where
$\{q^i\}$ denotes the residue   modulo  $r_j$.

\bf Let $H_j=r_j$. \rm To prove the mixing we  apply two properties of the  spacers: for  
$n<r_j$ we have 
\vspace{2mm}

{\bf (1)} $ - r_j \leq S_j(i,n):=\sum_{k=1}^{n} s_j(i+k)  - nH_j\leq r_j , \ \  i=1, 2, \dots, r_j-n-1$;
\vspace{2mm}

{\bf (2)}  for $ i\in\{1, 2, \dots, r_j-n-1\}$ all values $ S_j(i,n)$ are different.
\vspace{2mm}

Since  $$ \{q_j^i\} -\{q_j^{i+n}\} =S_j(i,n)=  S_j(m,n)=\{q_j^m\} -\{q_j^{m+n}\}$$  
implies  $$q^i-q^{i+n}=q^m-q^{m+n}, \ q^i=q^m,  \ i=m, \eqno (mod \ r_j) $$ we get  the  injectivity property (2).
\vspace{2mm}

\bf LEMMA. \it Let $r(j)>\delta r_j$ for a fixed  real $\delta\in (0,1)$, and $r(j)+n=r_j$.
If  $T$ is weakly mixing, then (1), (2)  imply $$Q(j)=\frac{ 1}{r(j)}\sum_{i=1}^{r(j)}T^{S_j(i,n)}\approx \Theta.$$
\rm

\bf COROLLARY. \it If a weakly mixing construction  satisfies (1),(2), then it is mixing.\rm
\vspace{2mm}

Proof.  We have to show that  for any $f\in L_2^0$  one has $\|Q(j)f\|\to 0$.  Otherwise there is a sequence $j_k$
such that $$Q({j_k})\to_w Q\neq \Theta. \eqno (\ast)$$
Defining a  measure $\eta$  by the formula
 $$\eta(A\times B)=\langle Q\chi_A , Q\chi_B\rangle_{L_2}$$
we see that $\eta=(T\otimes T)\eta$ and $\eta<<\mu\times\mu$.   The latter follows from 
 $$\mu(A)\mu(B)\gets \langle Q'(j)\chi_A , Q'(j)\chi_B\rangle \geq c \langle Q(j)\chi_A , Q(j)\chi_B\rangle\to c\eta(A\times B),$$
where $$Q'(j)=\frac{ 1}{2r_j}\sum_{i=-r_j}^{r_j}T^i,  \ \ c\geq \delta^2/4$$
(we remark that $|\{S_j(i,n): 1\leq i\leq r(j) \}|=r(j)\geq \delta r_j$, and  $-r_j\leq S_j(i,n)\leq r_j$).
Since  $T$ is weakly mixing,   $T\otimes T$ is ergodic, so $\eta=\mu\times\mu$, $Q=\Theta$. This contradicts $(\ast)$ and shows that $\|Q(j)f\|\to 0$.

Replacing stochastic spacers by algebraic ones and providing the weakly mixing property (i. e. the absence of eigenfunctions) we  get  mixing constructions. 

\bf Are algebraic constructions weakly mixing? \rm  
 Note that  the density  
of $\{S_j(i,1): 1\leq i\leq r_j \}$ in $[-r_j, r_j]$ is approx $0.5$. 
Taking this into account we can  see that  only the eigenvalue $-1$  could appear.  There are several simple "ergodic" ways to conserve (1),(2) and eliminate an  eigenfunction by   adding a little spacer.   We must avoid a situation in which for  most   $i$  one has  $(-1)^{S(i,1)}=(-1)^{h_j}$, hence,  we should be out of the same parity for most of $\{q^i\}-\{q^{i+1}\}$. With a pleasure and our thanks to Oleg  German  we present  his remark on  ``a parity of  the parities"  for  certain  pure constructions
that provides them the  weakly mixing property.


\bf O.N. German's arguments. \rm In 1962 Burgess proved that the minimal primitive root modulo a prime $r$ is 
$O(r^{0.25 +\eps})$,
where $\eps$  is a fixed arbitrarily small positive real number.  Hence for each $r_j$ large enough
we may choose $q_j$ to be less than $\sqrt{r_j}$. So, let $r$ be a large prime number (particularly, $r$ is
supposed to be odd), and let $q$ be a primitive root modulo $r$, such that $1 < q < \sqrt{r}$.  Let us
split the interval $[0, r)$ into the union of $q$ intervals of length $r/q$ and denote these intervals
as follows
$$I_k=[kr/q,(k+1)r/q, k=0,1,\dots,q-1.$$
Let us also denote by $\sigma_i$ the parity of the difference $\{q^i\}-\{q^{i+1}\}$, i.e.  $$\sigma_i=\{q^i\}-\{q^{i+1}\} mod\ 2,$$
and divide the set $M = \{1,2,\dots, r-1\}$ into the the two subsets
$$ M_0=\{i\in M: \sigma_i= 0\}, \ \
M_1=\{i\in M: \sigma_i= 1\}.$$
Suppose$q$ is odd. Then the parities of $\{q^{i}\}$ and $\{q^{i+1}\}$ coincide (i.e. $\sigma_i = 0$) if and only
if $\{q^{i}\}\in  I_k$ with even $k$. In this case we have
$$\left| |M_0|-|M_1|\right|\leq r/q +q-1<r/2+\sqrt{r}.$$
Here we have made use of the fact that the numbers of integer points contained in two
intervals of equal length differ at most by 1.
Suppose q is even. Then $\sigma_i = 0$ if and only if $\{q^{i}\}\in  I_k$ and the parity of $\{q^{i}\}$  coincides
with that of $k$. Hence in this case we have
$$\left| |M_0|-|M_1|\right|\leq q-1<\sqrt{r},$$
for in each $I_k$ the numbers of even points and odd points differ at most by 1.
Thus, if $r$ is large, then in case of an odd $q$, both $M_0$ and $M_1$ contain almost $1/4$ of all
the residues, and in case of an even $q$ those portions make almost $1/2$.


 E-mail: vryz@mail.ru
\end{document}